\documentclass[preprint,12pt]{elsarticle}

\usepackage{amssymb}
\usepackage{amscd,amsmath,amsfonts,graphicx}
\usepackage{amsthm}

\newcommand{\R}{\mathbb{R}}

\def\bgmat{\begin{displaymath}}
\def\endmat{\end{displaymath}}
\def\bgc{\begin{center}}
\def\endc{\end{center}}

\def\and{\mbox{ and }}

\def\bgeq{\begin{equation}}
\def\endeq{\end{equation}}
\begin{document}

\baselineskip=22pt \centerline{\Large \bf An Identity for Two
Integral Transforms}  \centerline{\Large \bf  Applied to the
Uniqueness of a Distribution} \centerline{\Large \bf
 via its Laplace--Stieltjes Transform}
\vspace{1cm} \centerline{Gwo Dong Lin $^{a,b}$ and\ Xiaoling Dou
$^c$} \centerline{$^a$ Hwa-Kang Xing-Ye Foundation, Taiwan}
\centerline{$^b$ Academia Sinica, Taiwan} \centerline{$^c$ Waseda
University, Japan ~} \vspace{0.5cm} \noindent {\bf Abstract.} It is
well known that the Laplace--Stieltjes transform of a nonnegative
random variable (or random vector) uniquely determines its
distribution function. We extend this uniqueness theorem by using
the M\"untz--Sz\'asz Theorem and
 the identity for the Laplace--Stieltjes and
Laplace--Carson transforms of a distribution function.  The latter
 appears for the first
time to the best of our knowledge.  In particular, if $X$ and $Y$
are two nonnegative random variables with joint distribution $H$,
then $H$ can be characterized by a suitable set of {\it countably
many values} of its bivariate Laplace--Stieltjes transform.   The
general high-dimensional case is also investigated. Besides, Lerch's
uniqueness theorem for conventional Laplace transforms is extended
as well. The  identity can be used to simplify the calculation of
{Laplace--Stieltjes transforms} when the underlying distributions
have {singular} parts. Finally, some examples are given to
illustrate the characterization results via the uniqueness theorem.

\vspace{0.5cm} \hrule \medskip \vspace{0.2cm}\noindent{\bf MSC2020
Mathematics Subject Classifications}:  62E10, 60E05, 46F12, 30E05.\\
\noindent{\bf Key words and phrases:} Lerch's Theorem,
M\"untz--Sz\'asz Theorem, Laplace transform, Laplace--Stieltjes
transform, Laplace--Carson transform, Characterization of distribution.\\
{\bf Postal addresses:}\\
Gwo Dong Lin, (1) Social and Data Science Research Center, Hwa-Kang
Xing-Ye Foundation, Taipei, Taiwan, and (2) Institute of Statistical
Science, Academia Sinica, Taipei 11529, Taiwan (E-mail:
 gdlin@stat.sinica.edu.tw)\\
 Xiaoling Dou, Center for Data Science,
Waseda University, 1-6-1 Nishiwaseda, Shinjuku-ku, Tokyo 169-8051,
Japan (E-mail: dou@waseda.jp)

\newpage \noindent {\bf 1. Introduction}\\
\indent The Laplace--Stieltjes transform  of  nonnegative random
variables (or their distribution functions) plays a crucial role in
the areas of theoretical and applied probability, especially in
survival analysis and reliability theory. The most important
fundamental property of the transform might be the uniqueness
theorem; namely, any two nonnegative random variables (or random
vectors) sharing the same Laplace--Stieltjes transform should have
the same distribution function
(see, e.g., Farrell 1976, pp.\,16--18).\\
\indent Another related integral transform of functions  is the
conventional Laplace transform which  is a powerful tool in solving
systems of both ordinary and partial differential equations through
the uniqueness theorem -- Lerch's Theorem (see Theorem A below,
Aghili and Parsa Moghaddam 2011, or van der Pol and Bremmer 1955).
In this paper, we prove  {\it the identity for the
Laplace--Stieltjes and Laplace--Carson transforms of a distribution
function} (the latter is in terms of the conventional Laplace
transform; see Theorems E, F and 4).  We will apply this identity to
extend the previous uniqueness theorem for Laplace--Stieltjes
transforms with the help of M\"untz--Sz\'asz Theorem (for the
latter, see Theorems B and C below). Moreover, the identity can be
used to simplify the calculation of {Laplace--Stieltjes transforms},
especially, when the underlying distributions
have {singular} parts. \\
\indent In Section 2, we formally define the integral transforms in
question and provide the needed preliminary results. In Section 3,
the main results for one- and two-dimensional cases are proved
(Theorems 1--3).  In particular, let $X$ and $Y$ be two nonnegative
random variables with joint distribution $H$, and let
$\{m_i\}_{i=1}^{\infty}$ and $\{n_j\}_{j=1}^{\infty}$ be two
sequences of strictly increasing positive real numbers satisfying
the conditions: $\sum_{i=1}^{\infty}1/m_i=\infty,\
\sum_{j=1}^{\infty}1/n_j=\infty.$ Then the joint distribution $H$ is
uniquely determined by the countably many values of bivariate
Laplace--Stieltjes transform: $\{{\bf E}[\exp(-m_iX-n_jY)]: i,
j=1,2,\ldots\}.$ The general high-dimensional case is investigated
in Section 4. Section 5 shows the advantage of the identity when
calculating {Laplace--Stieltjes transforms} of distributions with
{singular} parts. We give in Section 6 some examples to illustrate
the characterization results via the uniqueness theorem. Appendix
provides tedious calculations for the proof of the identity.
\medskip\\
\noindent {\bf 2. Preliminaries}\\
\indent  For a real-valued measurable function $f$ on ${\mathbb R}_+
:= [0,\infty)$, we denote by $L_f$ the (conventional) Laplace
transform of $f$:
\begin{eqnarray}L_f(s)=\int_0^{\infty}f(x)e^{-sx}{\rm d}x,\ s>0,
\end{eqnarray}
provided the integral exists.  The following uniqueness theorem is due to Lerch (1903);
see, e.g., Cohen (2007), p.\,23.\medskip\\
\noindent{\bf Theorem A} (Lerch's Theorem). {\it Suppose $f$ and $g$
are continuous on ${\mathbb R}_+$ and of exponential type $a\ge 0$
at infinity, namely, $|f(x)|={\cal O}(e^{ax})$ and
 $|g(x)|={\cal O}(e^{ax})$ as $x\to\infty.$ If $L_f(s)=L_g(s)$ for all $s>a$, then $f(x)=g(x)$ on ${\mathbb R}_+.$}\\
\indent In other words, the Laplace transform $L_f$ uniquely
determines the continuous function $f$ of  exponential type on
${\mathbb R}_+.$
 Note that to recover $f$ from $L_f$, we can apply the Post--Widder Inversion Formula:
 $$f(x)=\lim_{n\to \infty}\frac{(-1)^n}{n!}\left(\frac{n}{x}\right)^{n+1}L_f^{(n)}\!\left(\frac{n}{x}\right),\ x>0,$$
 where $L_f^{(n)}$ denotes the $n$th derivative of $L_f$ (see, e.g., Cohen 2007, p. 37).
 Conversely,  from Watson's lemma it follows the asymptotic equivalence
 \begin{eqnarray}L_f(s)\sim \sum_{n=0}^{\infty}\frac{f^{(n)}(0)}{s^{n+1}}~~\hbox {as}~~s\to \infty,
 \end{eqnarray}
 provided that $f\in C^{\infty}$ (infinitely continuously differentiable)
 in a neighborhood of zero  (see, e.g., Miller 2006, p.\,53). Interestingly,  if $f(x)=e^{-x}$,
 then we can replace the above approximation (2) by equality, namely,
 $L_f(s)=(1+s)^{-1}=\sum_{n=0}^{\infty}{f^{(n)}(0)}/{s^{n+1}}$ for all $s>1.$

One of our main purposes in this paper is to extend Lerch's Theorem
(Theorem A) to the case of measurable functions (see Theorem D
below). To do this, we need the following lemmas and
M\"untz--Sz\'asz Theorem (Theorems B and C; M\"untz 1914, Sz\'asz
1916). The latter extends the famous Weierstrass approximation
theorem in $C[0,1],$ the space of all continuous functions on
$[0,1]$ with supremum norm (see, e.g.,  P\'erez and Quintana 2008).
We provide detailed proofs here for completeness, although some of
them might be known with different approaches in the literature.

 \noindent{\bf Lemma A.} {\it Let
$f$ and $g$ be two Lebesgue integrable nonnegative functions on
$(0,\infty).$ Assume further that $\int_0^{\infty}f(x)e^{-sx}{\rm
d}x=\int_0^{\infty}g(x)e^{-sx}{\rm d}x$ for all
$s>0.$ Then $f(x)=g(x)$\, a.e.\,(almost everywhere) on $(0,\infty).$}\\
{\bf Proof.} Note that $\int_0^{\infty}f(x)\,{\rm
d}x=\int_0^{\infty}g(x)\,{\rm d}x =: A\in[0,\infty)$ by the
Dominated Convergence Theorem. If $A=0,$ then $f(x)=g(x)=0$ a.e. on
$(0,\infty),$ due to the assumptions $f,g\ge 0.$  Suppose now that
$A\in(0,\infty).$ Let us define two absolutely continuous
distributions:
\[F(x)=\frac{1}{A}\int_0^xf(t)\,{\rm d}t,\ \ x\ge 0,\quad\ G(y)=\frac{1}{A}\int_0^yg(t)\,{\rm d}t,\ \ y\ge 0.
\]
Moreover, let $X$ and $Y$ obey the distributions $F$ and $G,$
respectively, denoted $X\sim F$ and $Y\sim G.$ Then $X$ and $Y$ have
the same Laplace transforms by assumptions:
\[{\bf E}[\exp(-sX)]={\bf E}[\exp(-sY)],\ s\ge 0.\]
This means that $F=G$ and hence the difference between their
densities $f(x)/A-g(x)/A=0$\, a.e.\,on $(0,\infty).$ Equivalently,
$f(x)=g(x)$ a.e. on $(0,\infty).$ The proof is complete.

 \noindent{\bf Lemma B.} {\it Let $f$ be a
Lebesgue integrable function on $(0,\infty)$ and let
$\int_0^{\infty}f(x)e^{-sx}{\rm d}x=0$ for all $s> 0.$ Then
$f(x)=0$\, a.e. on
$(0,\infty).$}\\
{\bf Proof.} Define $f_+(x)=\max\{0, f(x)\}$ and
$f_-(x)=-\min\{0,f(x)\},\ \ x>0.$ Then $f=f_+-f_-$ and both $f_+$
and $f_-$ are nonnegative Lebesgue integrable functions on
$(0,\infty)$ satisfying
\[\int_0^{\infty}f_+(x)e^{-sx}{\rm d}x=\int_0^{\infty}f_-(x)e^{-sx}{\rm d}x,\ \ s>0.
\]
It then follows from Lemma A that $f_+(x)=f_-(x)$\, a.e. on
$(0,\infty)$ and hence $f(x)=f_+(x) - f_-(x)=0$\, a.e.  on
$(0,\infty).$ The proof is complete.

\noindent{\bf Lemma C.} {\it Let $f$ be a  measurable function on
$(0,\infty)$ such that $\int_0^{\infty}f(x)e^{-sx}{\rm d}x=0,\ \ s>
a$ for
some $a\ge 0.$ Then $f(x)=0$\, a.e. on $(0,\infty).$}\\
{\bf Proof.} Define $f_a(x)=f(x)e^{-ax},\ x> 0.$ Then $f_a$ is a
Lebesgue integrable function on $(0,\infty)$ satisfying
\[\int_0^{\infty}f_a(x)e^{-sx}{\rm d}x=0,\ \ s> 0.
\]
By Lemma B, we have $f_a(x)=0$\, a.e. on $(0,\infty)$ and so is $f.$
This completes the proof.

\noindent{\bf Theorem B} {(M\"untz--Sz\'asz Theorem in $C[0,1]$).}
{\it Suppose that $\Lambda=\{\lambda_k\}_{k=1}^{\infty}$ is a
sequence of positive and distinct real numbers satisfying
$\inf\Lambda>0$
 and
$\sum_{k=1}^{\infty}1/\lambda_k=\infty.$ Then the
span$\{1,x^{\lambda_1},x^{\lambda_2},\ldots\}$ is dense in $C[0,1].$
Equivalently, the collection of all finite linear combinations of
the functions $\{1,x^{\lambda_1},x^{\lambda_2},\ldots\}$ is dense in
$C[0,1],$ or, we also say that the set
$\{1,x^{\lambda_1},x^{\lambda_2},\ldots\}$ is total in $C[0,1].$}\\
{\bf Proof.} Case (I): $\sup\Lambda=\infty.$ Without loss of
generality  we may assume that $\lim_{k\to\infty}\lambda_k=\infty.$
We will elaborate von Golitschek's (1983) elementary constructive
proof. Let $q>0$ be a real number and $q\notin \Lambda.$ Define a
sequence of functions $\{Q_n\}$  by $q$ and $\{\lambda_k\}$:
$Q_0(x)=x^q,\ x\in[0,1],$ and for $n\ge 1,$
\[Q_n(x)=(\lambda_n-q)x^{\lambda_n}\int_x^1Q_{n-1}(t)t^{-(1+\lambda_n)}{\rm d}t,\ \ x\in[0,1].
\]
Then we have explicitly $Q_1(x)=x^q-x^{\lambda_1} =:
x^q-a_{1,1}x^{\lambda_1},\ x\in[0,1],$ and
\[Q_2(x)=x^q-\frac{\lambda_2-q}{\lambda_2-\lambda_1}x^{\lambda_1}-
\bigg(1-\frac{\lambda_2-q}{\lambda_2-\lambda_1}\bigg)x^{\lambda_2} =:
x^q-a_{1,2}x^{\lambda_1}-a_{2,2}x^{\lambda_2},\ x\in[0,1],\] where
$a_{2,2}=1-a_{1,2}.$  In general, for $n\ge 2,$ suppose $Q_{n-1}(x)$
is of the form:
\[Q_{n-1}(x)=x^q-\sum_{k=1}^{n-1}a_{k,n-1}x^{\lambda_k} =: x^q-P_{n-1}(x),\ \ x\in[0,1],\]
where $P_{n-1}(x)$ is a finite linear combination of
$x^{\lambda_1},x^{\lambda_2},\ldots,x^{\lambda_{n-1}}.$ Then we
carry out $Q_n(x)$ by definition and obtain
\begin{eqnarray*}
Q_{n}(x)&=&(\lambda_n-q)x^{\lambda_n}\int_x^1
\bigg(x^q-\sum_{k=1}^{n-1}a_{k,n-1}x^{\lambda_k}\bigg)t^{-(1+\lambda_n)}{\rm d}t\\
&=&x^q-\sum_{k=1}^{n}a_{k,n}x^{\lambda_k}=: x^q-P_n(x),\ \
x\in[0,1],
\end{eqnarray*}
where $P_{n}(x)$ is a finite linear combination of
$x^{\lambda_1},x^{\lambda_2},\ldots,x^{\lambda_{n}}$ with
coefficients
\[a_{k,n}=a_{k,n-1}\frac{\lambda_n-q}{\lambda_n-\lambda_k},\ k=1,2,\ldots,n-1,\
\quad a_{n,n}=1-\sum_{k=1}^{n-1}{a_{k,n}}.\] Hence, the
coefficients $\{a_{k,n}\}$ in $Q_n(x)$ can be derived from
$\{a_{k,n-1}\}$ iteratively. We now estimate the supremum of
$|Q_n(x)|$ on $[0,1]$, denoted $\|Q_n\|.$ By the definition of
$Q_n,$
\[\|Q_0\|=1,\quad \|Q_1\|\le \Big|1-\frac{q}{\lambda_1}\Big|,\]
\begin{eqnarray*} \|Q_n\|&\le&
\Big|1-\frac{q}{\lambda_n}\Big|\|Q_{n-1}\|\sup_{0\le x \le
1}|x^{\lambda_n}[1-x^{-\lambda_n}]|\\
&\le&
\Big|1-\frac{q}{\lambda_n}\Big|\|Q_{n-1}\|\le\prod_{k=1}^n\Big|1-\frac{q}{\lambda_k}\Big|\
 \longrightarrow 0\ \ \hbox{as}\ \ n\to\infty.
\end{eqnarray*}
The last conclusion is due to the assumptions that
$\lim_{k\to\infty}\lambda_k=\infty$ and
$\sum_{k=1}^{\infty}1/\lambda_k=\infty$ (Apostol 1975, p.\,209).
Namely, for any real $q>0,$ the monomials $x^q$ belong to the
closure of the span$\{x^{\lambda_k}: k=1,2,\ldots\}.$  Therefore,
the span$\{1,x^{\lambda_1},x^{\lambda_2},\ldots\}$ is dense in
$C[0,1]$ by the Weierstrass approximation theorem. The latter
asserts that every function $f\in C[0,1]$ is a uniform limit of
polynomials (Weierstrass 1885;  P\'erez and Quintana 2008).

 \noindent Case (II): $\sup\Lambda<\infty.$ Without loss of generality
 we may assume $\lim_{k\to\infty}\lambda_k=\lambda_*\in(0,\infty).$
Recall that $C[0,1]$ is dense in $L(0,1),$ the space of all Lebesgue
integrable functions on $(0,1).$ Suppose, on the contrary, that the
span$\{1,x^{\lambda_1},x^{\lambda_2},\ldots\}$ is not dense in
$C[0,1],$ and hence not dense in $L(0,1)$ either. Then by the
Hahn--Banach Theorem there exists a bounded nonzero measurable
function $g$ such that$\int_0^1x^{\lambda_k}g(x)\,{\rm d}x=0$ for
all $k.$ Define the complex-valued function
\[h(z)=\int_0^1x^{z}g(x)\,{\rm d}x,\ z\in\Pi=\{z: \hbox{Re}\ z>0\}.\]
Then $h$ is a bounded analytic function on the right half-plane
$\Pi$ and has zero value at the points $\lambda_k$ and the limit
$\lambda_*\in\Pi.$ This implies that $h(z)=0$ on $\Pi.$ By the
uniqueness theorem for Mellin transforms (or using Lemma C by
changing variables), $g(x)=0$ a.e. on $(0,1),$ a contradiction.
Therefore, the span$\{1,x^{\lambda_1},x^{\lambda_2},\ldots\}$ is
dense in $C[0,1].$
 The proof is complete.\\
 \indent We need some more notations.  Let $0\le a<b\le\infty$ and denote by $L(a,b)$ the space of all
 Lebesgue integrable functions on $(a,b)$.
 We say that the set of functions $\{f_n\}_{n=1}^{\infty}$ is complete in the space
 $L(a,b)$ if for any function $g\in L(a,b)$, the equalities
 $$\int_a^bf_n(x)g(x)\,{\rm d}x=0,\ \ n\in{\mathbb N} :=\{1,2,\ldots\},$$
together imply that $g(x)=0$\, a.e. on $(a,b)$ (Boas 1954, p.\,234).

\indent We have the following ramification of M\"untz--Sz\'asz
Theorem.

\noindent{\bf Theorem C.} {\it Let
$\Lambda=\{\lambda_k\}_{k=1}^{\infty}$ be a sequence of positive and
distinct real numbers satisfying $\inf\Lambda>0$
 and
$\sum_{k=1}^{\infty}1/\lambda_k=\infty.$ Then the set of functions
$\{x^{\lambda_k}\}_{k=1}^{\infty}$ is complete in
$L(0,1).$}\\
{\bf Proof.} Suppose $g\in L(0,1)$ and for $\lambda_k,$ we have the
equalities
\[\int_0^1t^{\lambda_k}g(t)\,{\rm d}t=0,\ \ k\in{\mathbb N}.
\]
Then we want to prove that $g(t)=0$ a.e. on $(0,1).$ By changing
variables $x=-\log t,$ we rewrite the above equalities in the form
\[0=\int_0^{\infty}e^{-\lambda_k x}g(e^{-x})e^{-x}{\rm d}x =: \int_0^{\infty}
e^{-\lambda_k x}h(x)\,{\rm d}x = L_h(\lambda_k),\ \ k\in{\mathbb N},
\]
where $h(x)=g(e^{-x})e^{-x}\in L(0,\infty)$ with Laplace transform
$L_h(\lambda)=\int_0^{\infty}e^{-\lambda x}h(x)\,{\rm d}x,
\lambda>0.$ Define the complex-valued function
\[L_h(z)= \int_0^{\infty}e^{-z x}h(x)\,{\rm d}x,\ \ z\in\Pi=\{z: \hbox{Re}\ z>0\}.
\]
Then $L_h$ is bounded and analytic on the right half-plane $\Pi$ and
vanishes at $\lambda_k$ for all $k\in{\mathbb N}.$  Set the function
\[H(z)=L_h\bigg(\frac{1+z}{1-z}\bigg)\ \hbox{for}\ z\in U=\{z: |z|<1\}.
\]
Then $H$ is bounded and analytic in $U.$
 Letting $\alpha_k=(\lambda_k-1)/(\lambda_k+1),$ we
have that $H(\alpha_k)=0$ and
$\sum_{k=1}^{\infty}(1-|\alpha_k|)=\infty.$ Therefore, $H(z)=0$ on
$U$ (Rudin 1987, p.\,312), or, equivalently, $L_h(z)=0$ on $\Pi.$
Particularly,  $L_h(\lambda)=0$ on $(0,\infty).$ By the uniqueness
theorem for Laplace transforms of measurable functions (see Lemma C
above), we conclude that $h(x)=0$\, a.e. on $(0,\infty)$ and hence
$g(t)=0$\, a.e. on $(0,1).$ This completes the proof.

It is seen that the next theorem improves significantly both Theorem
A and Lemma C.

 \noindent {\bf Theorem D.} {\it Let  $f$ be a
measurable function on $(0,\infty)$ and let $\{n_j\}_{j=1}^{\infty}$
be a sequence of positive and distinct increasing real numbers
satisfying $\sum_{j=1}^{\infty}1/n_j=\infty.$ Assume further that
$\int_0^{\infty}f(x)e^{-n_j x}{\rm d}x=0$ for all $j\in{\mathbb N}.$
Then $f(x)=0$\, a.e. on $(0,\infty).$ }\\
{\bf Proof.} We rewrite, by changing variables $t=e^{-x}$,
$$\int_0^{\infty}f(x)e^{-n_j x}{\rm d}x=\int_0^{\infty}[f(x)e^{-n_1x}]e^{-(n_j-n_1) x}{\rm d}x=\int_0^1h(t)t^{n_j-n_1}{\rm d}t=0,\ j=2,3,\ldots,$$
where the function $h(t)=f(-\ln t)t^{n_1-1}\in L(0,1)$ since
$\int_0^{\infty}f(x)e^{-n_1x}{\rm d}x=0.$ Then by Theorem C,
$h(t)=0$\, a.e. on $(0,1)$ because
$\sum_{j=2}^{\infty}1/(n_j-n_1)\ge\sum_{j=2}^{\infty}1/n_j=\infty,$
and hence $f(x)=0$\, a.e. on $(0,\infty).$ The proof is complete.
\medskip\\
\indent For recent developments on the M\"untz--Sz\'asz Theorem, see
Erd\'elyi and Johnson (2001), Erd\'elyi (2005), Almira (2007) and
the references therein.

 \indent We now consider two nonnegative
random variables $X$ and $Y$ having joint distribution $H$ with
marginals $F$ and $G$, that is, $(X, Y)\sim H$, $X\sim F$ and $Y\sim
G.$ Denote by ${\cal L}_F$ and ${\cal L}_H$ the Laplace--Stieltjes
transforms of $X\sim F$ and $(X,Y)\sim H,$ respectively. Formally,
\begin{eqnarray}{\cal
L}_F(s)&=&{\bf E}[\exp(-sX)]=\int_{0}^{\infty}e^{-sx}{\rm d}F(x)=\int_{[0,\infty)}e^{-sx}{\rm d}F(x)\\
&=&\lim_{\varepsilon\downarrow 0}\int_{-\varepsilon}^{\infty}e^{-sx}{\rm d}F(x)=F(0)+\int_{(0,\infty)}e^{-sx}{\rm d}F(x),\ \ s>0,\nonumber\\
 {\cal
L}_H(s,t)&=&{\bf E}[\exp(-sX-tY)]=\int_0^{\infty}\!\!\int_0^{\infty}e^{-sx-ty}{\rm d}H(x,y),\ \
s,t>0.
\end{eqnarray}
 Also, analogously to  (1),
denote by $L_H$ the (conventional) Laplace transform of bivariate
$H$:
\begin{eqnarray}L_H(s,t)=\int_0^{\infty}\!\!\int_0^{\infty}H(x,y)e^{-sx-ty}{\rm d}x{\rm d}y,\ \
s,t>0.
\end{eqnarray} \indent It is known that ${\cal L}_F\in
C^{\infty}((0,\infty))$ is completely monotone and that for each
continuity point $x>0$ of $F$,
\begin{eqnarray}F(x)=\lim_{n\to\infty}\sum_{k\le
nx}(-1)^k\frac{n^k}{k!}{\cal L}_F^{(k)}(n)=
\lim_{n\to\infty}\sum_{k=0}^n\frac{(-1)^k}{k!}\left(\frac{n}{x}\right)^k{\cal
L}_F^{(k)}\!\left(\frac{n}{x}\right)
\end{eqnarray} (see, e.g., Feller 1971, pp.\,439--440). The second equality in (6) also follows from the facts: (i) Theorem E
below, (ii) the Post--Widder Inversion Formula and (iii) Leibniz's
rule. Besides, there are several interesting and useful
relationships between the Laplace--Stieltjes transform and the
conventional Laplace transform. Note first that if $F$ has a density
$f$, then ${\cal L}_F=L_f$ by definition, and that if the bivariate
$H$ has a joint density $h$, then ${\cal L}_H=L_h.$ The following
two identities can be derived by using the existing results and will
be used in the sequel.

\noindent{\bf Theorem E.} {\it For any distribution $F$ on ${\mathbb R}_+,\ $
\[{\cal L}_F(s)=sL_F(s)\ \ for\ all\ \ s>0.\]}
{\bf Proof.} For $s>0,$ using integration by parts we have
\begin{eqnarray}&~&{\cal L}_F(s)=\int_0^{\infty}e^{-sx}{\rm d}F(x)=F(0)+\int_{(0,{\infty})}e^{-sx}{\rm d}F(x)\nonumber\\
&=&F(0)+e^{-sx}F(x)|_{0^+}^{\infty}+s\int_{(0,{\infty})}F(x)e^{-sx}{\rm d}x=sL_F(s).
\end{eqnarray} Also, note that  the identity (7) has an equivalent form in terms of the survival function $\overline{F}:$
\[(1-{\cal L}_F(s))/s=\int_0^{\infty}e^{-sx}\overline{F}(x)\,{\rm d}x,\ \ s>0,
\]
where $\overline{F}(x)={\bf P}(X>x)=1-F(x),\ x\ge
0$ (see  Lin 1998, Lemma 1, or Feller 1971, p.\,435).\medskip\\
{\bf \bf Theorem F.} {\it For any bivariate distribution $H$ on
${\mathbb R}_+^2,\ $
\[ {\cal L}_H(s,t)=stL_H(s,t)\ \ for\ all\ \ s, t>0.\]}
{\bf Proof.} From Theorem E it follows that the identity in question
is equivalent to
\begin{eqnarray}{\cal L}_H(s,t)=st\int_{0}^{\infty
}\!\!\int_{0}^{\infty }\overline{H} (x,y)e^{-sx-ty}{\rm d}x{\rm d}y-1
 +{\cal L}_F(s)+{\cal L}_G(t),\ s,\,t\ge 0\end{eqnarray}
(by an expansion of the double integral on the RHS of (8)),  where
the joint survival function
\[\overline{H}(x,y)={\bf P}(X>x,\ Y>y)=1-F(x)-G(y)+H(x,y),\ x,y\ge 0.\]
 The identity (8) is exactly Corollary 2 in  Lin et al.\,(2016). The proof is complete.
\medskip\\
\indent The RHS ($sL_F(s)$) of the identity in Theorem E is called
the Laplace--Carson transform of a function $F$ (compare (1) and
(3)). For definition of the latter transform, see, e.g., Carson
(1919), Ditkin and Prudnikov (1962) and Donolato (2002). In other
words, Theorems E and F claim {\it the  identity for the
Laplace--Stieltjes and the Laplace--Carson transforms of a
distribution function} in the first two dimensions. To the best of
our knowledge, at least the bivariate identity in Theorem F appears
for the first
time.\medskip\\
 \noindent {\bf 3. Main results: one- and two-dimensional
cases}

\indent  We restate Theorem D as follows.

\noindent {\bf Theorem 1.} {\it Let  $f, g$ be two measurable
functions on $(0,\infty)$ and let $\{n_j\}_{j=1}^{\infty}$ be a
sequence of positive and distinct increasing real numbers satisfying
$\sum_{j=1}^{\infty}1/n_j=\infty.$ If $\int_0^{\infty}f(x)e^{-n_j
x}{\rm d}x=\int_0^{\infty}g(x)e^{-n_j x}{\rm d}x$ $($finite$)$  for
all $j\in{\mathbb N},$ then $f(x)=g(x)$\ a.e. on $(0,\infty).$}

\noindent{\bf Theorem 2.} {\it Let $0\le X\sim F$ and let the real
numbers  $0<m_1<m_2<\cdots$ satisfy
$\sum_{i=1}^{\infty}1/m_i=\infty.$ Then the distribution $F$ is
uniquely determined by the countable set of
values $\{{\cal L}_F(m_i)\}_{i=1}^{\infty}$ of its Laplace--Stieltjes transform.} \\
{\bf Proof.} Suppose $0\le X_1\sim F_1$, $0\le X_2\sim F_2$ and
${\cal L}_{F_1}(m_i)={\cal L}_{F_2}(m_i)$, $i\in{\mathbb N}.$ Then
we want to show that $F_1=F_2$ under the condition
$\sum_{i=1}^{\infty}1/m_i=\infty.$ By Theorem E and the assumptions,
we have the equalities:
$$\int_0^{\infty}F_1(x)e^{-m_ix}{\rm d}x=\int_0^{\infty}F_2(x)e^{-m_ix}{\rm d}x\ \ (\hbox{finite})\ \ \forall\  i\in{\mathbb N}.$$
It then follows from Theorem 1 that $F_1(x)=F_2(x)$\ a.e. on
$[0,\infty).$
 and hence $F_1=F_2$ due to the right
continuity  of distributions. The proof is complete. \\
\indent An alternative proof of Theorem 2 was given in Lin (1993,
Lemma 4), in which two other sufficient conditions were provided:\\
(i) $\lim_{i\to\infty}m_i=m_0\in(0,\infty),$ and \\(ii)
$\lim_{i\to\infty}m_i=0,\ \sum_{i=1}^{\infty}m_i=\infty$.\\ These
three  conditions also apply to the high-dimensional cases, but for
simplicity, we consider only the strictly monotone  sequence below.
\medskip\\
 \noindent{\bf Theorem 3.} {\it Let $0\le X\sim F,\ 0\le Y\sim G$ and
$(X,Y)\sim H.$ Assume further that the two sequences
$\{m_i\}_{i=1}^{\infty}$ and $\{n_j\}_{j=1}^{\infty}$  of real numbers satisfy\\
 (i)
$0<m_1<m_2<\cdots$ with $\sum_{i=1}^{\infty}1/m_i=\infty,$ and\\
 (ii)
$0<n_1<n_2<\cdots$ with $\sum_{j=1}^{\infty}1/n_j=\infty.$\\ Then
the bivariate distribution $H$ is uniquely determined by the
countably many
values  of its Laplace--Stieltjes transform: $\{{\cal L}_H(m_i,n_j): i,j=1,2,\ldots\}.$}\\
{\bf Proof.} For $k=1,2,$ suppose $0\le X_k\sim F_k,\ 0\le Y_k\sim
G_k,$  $(X_k,Y_k)\sim H_k,$ and ${\cal L}_{H_1}(m_i,n_j)={\cal
L}_{H_2}(m_i,n_j),$ $i,j\in{\mathbb N}.$ Then we want to prove that
$H_1=H_2$ under the conditions $\sum_{i=1}^{\infty}1/m_i=\infty$ and
$\sum_{j=1}^{\infty}1/n_j=\infty.$ By Theorem F and the assumptions,
we have the equalities:
\begin{eqnarray}(m_in_j)\int_0^{\infty}\left[\int_0^{\infty}H^*(x,y)e^{-m_ix}{\rm d}x\right]e^{-n_jy}{\rm d}y=0,\
i,j\in{\mathbb N},
\end{eqnarray} where the function $H^*=H_1-H_2.$  For fixed $y,$ let $H^*_y$ denote a function $H^*(x,y)$ of $x.$
It then follows from (9) and Theorem D that
\begin{eqnarray}{L}_{H_y^*}(m_i) :=
\int_0^{\infty}H^*(x,y)e^{-m_ix}{\rm d}x=0\ \  \forall\ y\ge 0,\
i\in{\mathbb N},\end{eqnarray} because
$\sum_{j=1}^{\infty}1/n_j=\infty$ and $H^*(x,y)$ is right continuous
in $y.$ By (10) and Theorem D again, we have that $H^*(x,y)=0\ \
\forall\ x, y\ge 0,$ due to the assumption
$\sum_{i=1}^{\infty}1/m_i=\infty$ and the right continuity of
$H^*(x,y)$ in $x.$ Therefore, $H_1=H_2$. This completes the proof.
\medskip\\
\indent  It is seen that the identity for two integral transforms
plays a crucial role in the proofs of Theorems 2 and 3. Motivated by
the results of the one- and two-dimensional cases, we next consider
the $n$-dimensional case with $n\ge 3,$ which is much more
complicated.
\medskip\\
\noindent {\bf 4. The general result: $n$-dimensional case with $n\ge 3$}\\
\indent Consider the random variables $0\le X_i\sim H_i,\
i=1,2,\ldots, n\ (\ge 3),$ and suppose that they have the joint
distribution $H$ on $\mathbb R_+^n.$ For $s_i>0,\ i=1,2,\ldots,n,$
denote the Laplace--Stieltjes transform of $H$ by \[{\cal
L}_H(s_1,s_2,\ldots,s_n)={\bf
E}\bigg[\exp\Big(-\sum_{i=1}^ns_iX_i\Big)\bigg]\]
\[=\int_0^{\infty}\!\!\int_0^{\infty}\!\!\cdots\int_0^{\infty}\exp\Big({-\sum_{i=1}^ns_ix_i}\Big){\rm d}H(x_1,x_2,\ldots,x_n),\]
and the (conventional) Laplace transform of $H$ by
\[{L}_H(s_1,s_2,\ldots,s_n)
=\int_0^{\infty}\!\!\int_0^{\infty}\!\!\cdots\int_0^{\infty}
H(x_1,x_2,\ldots,x_n)\exp\Big({-\sum_{i=1}^ns_ix_i}\Big){\rm d}x_1{\rm d}x_2\cdots {\rm d}x_n.
\]
Under the above setting, we extend Theorems E and F to the following identity in dimension $n,$ where $n\ge 3.$\\
{\bf Theorem 4.} {\it The Laplace--Stieltjes transform of the joint
distribution $H$ on $\mathbb{R}_+^n$ satisfies
 \begin{eqnarray} {\cal
L}_H(s_1,s_2,\ldots,s_n)=\bigg(\prod_{i=1}^ns_i\bigg)L_H(s_1,s_2,\ldots,s_n),\
\ s_i>0,\ i=1,2,\ldots,n.\end{eqnarray}} {\bf Proof.} (I)  We first
consider the special case:  $H$ is absolutely continuous with
positive marginal densities on $(0,\infty),$  because its proof is
simpler and easier to understand than that in the general case.  In
this case, we prove the identity (11) by induction on $n.$ It
follows from Theorems E and F that the identity (11) holds true for
$n=1$ and $n=2.$ Now, suppose it holds true for $n=m\ge 2,$ then we
want to prove the identity for $n=m+1.$ Denote the density of $H$ by
$h$ and its $j$th marginal density by $h_j$. Then the marginal
density of the first $m+1$ components $X_1,X_2,\ldots, X_m,X_{m+1}$
is written as follows: \[h(x_1,x_2,\ldots,x_{m},x_{m+1})
=h(x_1,x_2,\ldots,x_{m}|x_{m+1})h_{m+1}(x_{m+1}),\] and by the
absolutely continuous condition on $H,$
\begin{eqnarray*}
& &{\cal
L}_H(s_1,s_2,\ldots,s_{m+1})=L_h(s_1,s_2,\ldots,s_{m+1})\\
&=&\int_0^{\infty}\!\!\int_0^{\infty}\!\!\cdots\int_0^{\infty}\exp\Big({-\sum_{i=1}^{m+1}s_ix_i}\Big)
h(x_1,x_2,\ldots,x_m|x_{m+1})h_{m+1}(x_{m+1})\,{\rm d}x_1{\rm d}x_2\cdots {\rm d}x_{m+1}\\
&=&\int_0^{\infty}\exp\big({-s_{m+1}x_{m+1}}\big)h_{m+1}(x_{m+1})\\
& & \times\left[\int_0^{\infty}\!\!\cdots\int_0^{\infty}\exp\Big({-\sum_{i=1}^{m}s_ix_i}\Big)
h(x_1,x_2,\ldots,x_m|x_{m+1})\,{\rm d}x_1{\rm d}x_2\cdots {\rm d}x_m\right]{\rm d}x_{m+1}.
\end{eqnarray*}
By the assumption on the case $n=m$, the factor in square brackets
is equal to
$$\bigg(\prod_{i=1}^ms_i\bigg)\int_0^{\infty}\!\!\cdots\int_0^{\infty}
H(x_1,x_2,\ldots,x_m|X_{m+1}=x_{m+1})\exp\Big({-\sum_{i=1}^{m}s_ix_i}\Big)\,{\rm d}x_1{\rm d}x_2\cdots
{\rm d}x_m.$$ Next, rewrite
\begin{eqnarray*}
& &H(x_1,x_2,\ldots,x_m|X_{m+1}=x_{m+1})h_{m+1}(x_{m+1})\\
&=&h_{m+1}(x_{m+1}|X_1\le
x_1,\ldots,X_m\le x_m){\bf P}(X_1\le
x_1,\ldots,X_m\le x_m).
\end{eqnarray*}
Then  we have, by Theorem E,
\begin{eqnarray*}&
&\int_0^{\infty}e^{-s_{m+1}x_{m+1}}h_{m+1}(x_{m+1}|X_1\le
x_1,\ldots,X_m\le x_m){\rm d}x_{m+1}\\&=&s_{m+1}\int_0^{\infty}{\bf P}(X_{m+1}\le
x_{m+1}|X_1\le x_1,\ldots,X_m\le x_m)e^{-s_{m+1}x_{m+1}}\,{\rm d}x_{m+1},
\end{eqnarray*}
and hence the required conclusion follows by Fubini's theorem,
because
\begin{eqnarray*}&~&{\bf P}(X_{m+1}\le x_{m+1}|X_1\le x_1,\ldots,X_m\le x_m){\bf P}(X_1\le
x_1,\ldots,X_m\le x_m)\\
&=&H(x_1,x_2,\ldots,x_{m+1}).\end{eqnarray*} (II) We
next treat the general case, in which no smoothness  conditions on
distributions are assumed. Again, we prove the identity (11) by
induction on $n.$ Suppose it holds true for $n=1,2,\ldots,m\ge 2,$
then we want to prove the identity for $n=m+1.$

For fixed $s_i>0,\ i=1,2,\ldots,m+1,$ define the function $K$ on
$\R_+^{m+1}$ by
\[K(x_1,x_2,\ldots,x_{m+1})=\prod_{i=1}^{m+1} \big(1-\exp(-s_ix_i)\big),\ x_i\ge 0,\
i=1,2,\ldots,{m+1}.\] Then the Lebesgue--Stieltjes integral
\begin{eqnarray*}
&~&{\bf E}[K(X_1,X_2,\ldots,X_{m+1})]\\
&=&\int_0^{\infty}\!\!\int_0^{\infty}\!\!\cdots\int_0^{\infty}K(x_1,x_2,\ldots,x_{m+1})\,{\rm d}H(x_1,x_2,\ldots,x_{m+1})\\
&=&\lim_{R\to\infty}\int_0^{R}\!\!\int_0^{R}\!\!\cdots\int_0^{R}K(x_1,x_2,\ldots,x_{m+1})\,{\rm d}{H}(x_1,x_2,\ldots,x_{m+1})\\
&=&(-1)^{m+1}\lim_{R\to\infty}\int_0^{R}\!\!\int_0^{R}\!\!\cdots\int_0^{R}K(x_1,x_2,\ldots,x_{m+1})\,{\rm d}\overline{H}(x_1,x_2,\ldots,x_{m+1}).
\end{eqnarray*}
Here the joint survival function is
\begin{eqnarray}&~&\overline{H}(x_1,x_2,\ldots,x_{m+1})={\bf P}(X_1>x_1,X_2>x_2,\ldots,X_{m+1}>x_{m+1})\nonumber\\
&=&1-\sum_{i=1}^{m+1}H_i(x_i)+\sum_{1\le i_1<i_2\le m+1}H_{i_1i_2}(x_{i_1},x_{i_2})-\sum_{1\le i_1<i_2<i_3\le m+1}H_{i_1i_2i_3}(x_{i_1},x_{i_2},x_{i_3})\nonumber\\
&~~&+\cdots+(-1)^{m+1}H(x_1,x_2,\ldots,x_{m+1}),
\end{eqnarray}
and the $k$-dimensional marginal distribution, that is, the joint
distribution of $X_{i_1}, X_{i_2},\ldots,X_{i_k},$ is
\begin{eqnarray}H_{i_1i_2\ldots
i_k}(x_{i_1},x_{i_2},\ldots,x_{i_k})={\bf P}(X_{i_1}\le
x_{i_1},X_{i_2}\le x_{i_2},\ldots,X_{i_k}\le x_{i_k}).\end{eqnarray}
Formula (12) for $\overline H$ is valid for any collection of random
variables. The proof is standard, it is based on the well-known
inclusion-exclusion representation for the union of arbitrary
collection of random events, then we take probability and follow
induction arguments. Details can be seen in many books in
Probability Theory; see, e.g., Ross (2014), p.\,6.


 Using
the multidimensional integration by parts (see Young 1917, Section
9, or Zaremba 1968, Proposition 2) and proceeding in a similar way
as  in Lin et al.\,(2016), pp.\,3--4, we have
\begin{eqnarray}
&~&(-1)^{m+1}\lim_{R\to\infty}\int_0^{R}\!\!\int_0^{R}\!\!\cdots\int_0^{R}K(x_1,x_2,\ldots,x_{m+1})\,{\rm d}\overline{H}(x_1,x_2,\ldots,x_{m+1})\nonumber\\
&=&(-1)^{2({m+1})}\lim_{R\to\infty}\int_0^{R}\!\!\int_0^{R}\!\!\cdots\int_0^{R}\overline{H}(x_1,x_2,\ldots,x_{m+1})\,{\rm d}K(x_1,x_2,\ldots,x_{m+1})\\
&=&\bigg(\prod_{i=1}^{m+1}s_i\bigg)\int_0^{\infty}\!\!\int_0^{\infty}\!\!\cdots\int_0^{\infty}
\overline{H}(x_1,x_2,\ldots,x_{m+1})\exp\Big({-\sum_{i=1}^{m+1}s_ix_i}\Big){\rm d}x_1{\rm d}x_2\cdots {\rm d}x_{m+1},\nonumber
\end{eqnarray}
in which to derive the identity (14),  we apply the facts that
$K(x_1,x_2,\ldots,x_{m+1})=0$ if $x_i=0$ for some $i$ and
that\[\overline{H}(x_1,x_2,\ldots,x_{m+1})\to 0\ \ \hbox{if}\ x_j\to
\infty\ \hbox{for some}\ j.\]  Therefore, for $s_i>0,\
i=1,2,\ldots,{m+1},$ the above Lebesgue--Stieltjes integral becomes
\begin{eqnarray}
&~&{\bf E}\bigg[\prod_{i=1}^{m+1} \big(1-\exp\big({-s_iX_i}\big)\big)\bigg]\nonumber\\
&=&\!\!\!
\bigg(\prod_{i=1}^{m+1}s_i\bigg)\int_0^{\infty}\!\!\int_0^{\infty}\!\!\cdots\int_0^{\infty}
\overline{H}(x_1,x_2,\ldots,x_{m+1})\exp\Big({-\sum_{i=1}^{m+1}s_ix_i}\Big){\rm d}x_1{\rm d}x_2\cdots
{\rm d}x_{m+1}.~~~
\end{eqnarray}
The last identity reduces to  the following one  by the assumption
on induction and by a tedious calculation (using (12) and canceling
the same terms on both sides of (15)): \begin{eqnarray} &~&{\bf
E}\Big[\exp\Big({-\sum_{i=1}^{m+1}s_iX_i}\Big)\Big]\nonumber\\
&=&\!\!\!\bigg(\prod_{i=1}^{m+1}s_i\bigg)\int_0^{\infty}\!\!\int_0^{\infty}\!\!\cdots\int_0^{\infty}
{H}(x_1,x_2,\ldots,x_{m+1})\exp\Big({-\sum_{i=1}^{m+1}s_ix_i}\Big){\rm
d}x_1{\rm d}x_2\cdots {\rm d}x_{m+1}.~~~
\end{eqnarray}
 (See
Appendix  for a detailed proof of the identity (16).) The proof is
complete.

 Using the crucial identity (11) and mimicking the proof of Theorem 3, we have, by induction,
 the following uniqueness theorem for Laplace--Stieltjes transforms in
 the high-dimensional case.  The proof is omitted.\medskip\\
{\bf Theorem 5.} {\it Let $\{m_{i,j}\}_{j=1}^{\infty},\
i=1,2,\ldots,n,$ be the $n\ (\ge 3)$ sequences of real numbers
satisfying
\[0<m_{i,1}<m_{i,2}<\cdots\ \ \hbox{and}\ \ \sum_{j=1}^{\infty}\frac{1}{m_{i,j}}=\infty\ \
\hbox{for\ all}\ \ i=1,2,\ldots,n.\] Then  any $n$-dimensional
distribution $H$ on $\R_+^n$ is uniquely determined by the countably
many values of its Laplace--Stieltjes transform: \[\{{\cal
L}_H(s_1,s_2,\ldots,s_n): s_i=m_{i,j},\
i=1,2,\ldots,n,\,j\in{\mathbb N}\}.\]} \indent The following is a
special case of Theorems 2, 3 and 5, because
$\sum_{j=1}^{\infty}1/j=\infty.$

\noindent{\bf Corollary 1.}  {\it Any $n$-dimensional distribution
$H$ on $\R_+^n$ is uniquely determined by the countably many values
of its Laplace--Stieltjes transform: $\{{\cal
L}_H(s_1,s_2,\ldots,s_n): s_i=j,\ i=1,2,\ldots,n,\,j\in{\mathbb
N}\}.$ In other words, the set of values ${\cal
L}_H(s_1,s_2,\ldots,s_n)$ at the lattice points in ${\mathbb N}^n$
characterizes the distribution $H.$}
\medskip\\
\noindent {\bf 5. Application to calculation of Laplace--Stieltjes
transforms}

In general, the calculation of the Laplace--Stieltjes transform
${\cal L}_H$ is much more complicated than that of the conventional
Laplace transform $L_H,$ especially, when the underlying
distribution $H$ has a {\it singular} part, which often happens in
survival analysis. The identity (11) then can be used to simplify
the calculation. To illustrate this advantage, let us consider the
bivariate-lack-of-memory (BLM) distribution defined below.

Let $X\sim F$ and $Y\sim G$ be two positive random variables having
joint distribution $H.$ Namely, $(X,Y)\sim H$ with marginals $F$ and
$G.$ We say that $(X,Y)$ has a BLM distribution if it satisfies the
BLM property:
\[{\bf P}(X>x+t,\,Y>y+t|\ X>t,\,Y>t)={\bf P}(X>x,\,Y>y),\
x,y,t\ge 0,\] or, equivalently, if its survival function $\overline{H}$ satisfies the functional equation:
\[{\overline{H}(x+t,\,y+t)=\overline{H}(x,y)\overline{H}(t,t),\
\forall\ x,y,t\ge0.} \] Explicitly,
 the survival
function $\overline{H}$  can be written in the form
\begin{equation*}
   { \overline{H}(x,y)=\left\{\begin{array}{cc}
                    e^{-\theta y}\,\overline{F}(x-y), & x\ge y\ge 0 \vspace{0.1cm}\\
                    e^{-\theta x}\,\overline{G}(y-x),  & y\ge x\ge 0,
                  \end{array}\right.}
\end{equation*}
where $\theta$  is a positive constant (see {Marshall and Olkin
1967, or Barlow and Proschan 1981, p.\,130).  For convenience, we
denote $H=BLM(F,G,\theta),$ which has a singular part on the line
$x=y$ with probability $p(\theta):=\big(f(0)+g(0)\big)/\theta-1\ge
0,$ where $f(0)=\lim_{\varepsilon\to
0^+}F(\varepsilon)/\varepsilon,$ and $g(0)=\lim_{\varepsilon\to
0^+}G(\varepsilon)/\varepsilon$ (see Remark 2 in Lin et al.\,2019).

When $p(\theta)>0,$ $H$ has a singular part and it is hard to
calculate directly the {Laplace--Stieltjes transform} ${\cal
L}_H(s,t)={\bf E}\left[\exp({-sX-tY})\right],\ s, t>0.$ Fortunately,
applying the helpful  identity, it suffices to carry out
\begin{eqnarray*}
stL_H(s,t)&=&st\int_0^{\infty}\int_0^{\infty}H(x,y)\exp({-sx-ty})\,{\rm d}x{\rm d}y\\
&=&st\int_0^{\infty}\int_0^{\infty}[\overline{H}(x,y)-1+F(x)+G(y)]\exp({-sx-ty})\,{\rm d}x{\rm d}y.
\end{eqnarray*}
The advantage is that we now can ignore the singular part of $H$ (or
$\overline{H}$), because its two-dimensional Lebesgue measure is
zero. The final result is
\begin{eqnarray*}{\cal L}_H(s,t)=\frac{1}{\theta+s+t}\left[(\theta+s){\cal L}_F(s)+(\theta+t){\cal L}_G(t)\right]-\frac{\theta}{\theta+s+t},\
s, t>0.
\end{eqnarray*}
(See Lin et al.\,2019, Theorem 2, for a different proof.)
\medskip\\
\noindent {\bf 6. Application to the characterization of distributions}\\
\indent To illustrate the use of Theorem 5, we will consider
 some frequently used
Laplace--Stieltjes transforms below. Denote by $p_j$ the $j$th prime
number (we have $p_1=2,\ p_2=3,\ p_3=5,\ldots$) and denote by ${\cal
P} :=\{p_j\}_{j=1}^{\infty}$ the sequence of all prime numbers. Then
it is known that $\sum_{j=1}^{\infty} 1/p_j=\infty$ because $p_j\sim
j\ln j$ as $j\to \infty$ (see, e.g., Apostol 1976, p.\,80).

Now we present eight examples in which $F$ is a $1$-dimensional
distribution and $H$ is a $2$-dimensional distribution in Examples
3--7, while $3$-dimensional distribution in Example 8. More examples
can be found in Balakrishnan and Lai (2009). Notice that based on
our results, we are in a position to formulate each example as a
direct statement involving a specific distribution.

\noindent {\bf Example 1.} The Laplace--Stieltjes transform ${\cal
L}_F$ satisfies \[{\cal L}_F(p_j)=\frac{\lambda}{\lambda+p_j},\
j\in{\mathbb N},\ {\rm where}\ \ \lambda>0,\] if and only if $F$ is
the exponential distribution with mean $1/\lambda.$ More generally,
\[{\cal L}_F(p_j)=\Big[\frac{\lambda}{\lambda+p_j}\Big]^q,\ j\in{\mathbb N},\ {\rm where}\ \
\lambda, q>0,\] if and only if $F$ is the Gamma distribution with
mean
$q/\lambda$ and shape parameter  $q$, namely, $F$ has a density $f(x)=[\lambda^q/\Gamma(q)]x^{q-1}\exp(-\lambda x),\ x\ge 0.$ \medskip\\
{\bf Example 2.}  The Laplace--Stieltjes transform ${\cal L}_F$
satisfies \[{\cal L}_F(p_j)=\exp[-p_j^{\alpha}],\ j\in{\mathbb N},\
{\rm where}\ \ \alpha\in (0,1),\] if and only if $F$ is the positive
stable distribution with density function
$$f_{\alpha}(x)=-\frac{1}{\pi x}\sum_{k=1}^{\infty}\frac{\Gamma(\alpha k+1)}{k!}
\left(-x^{-\alpha}\right)^k\sin(\alpha k\pi),\ \ x>0$$
 (see Feller 1971, p.\,583, and Hougaard 1986).\medskip\\
{\bf Example 3.} The Laplace--Stieltjes transform ${\cal L}_H$
satisfies
$${\cal
L}_H(p_i,p_j)=\frac{(\lambda+p_i+p_j)(\lambda_1+\lambda_{12})(\lambda_2+\lambda_{12})+p_ip_j\lambda_{12}}
{(\lambda+p_i+p_j)(\lambda_1+\lambda_{12}+p_i)(\lambda_2+\lambda_{12}+p_j)},\
\ i,j\in{\mathbb N},$$ where $\lambda_1,\lambda_2, \lambda_{12}$ are
positive constants and $\lambda=\lambda_1+\lambda_2+\lambda_{12},$
if and only if $H$ is the Marshall--Olkin bivariate exponential
(BVE) distribution with survival function
${\overline{H}(x,y)}={\exp[-\lambda_1 x-\lambda_2
y-\lambda_{12}\max\{x,y\}]},\ x,y\ge 0$\  \ (see Marshall and Olkin
1967).\medskip\\
{\bf Example 4.} The Laplace--Stieltjes transform ${\cal L}_H$
satisfies
$${\cal L}_H(p_i,p_j)=\frac{1}{\alpha+\beta+p_i+p_j}\left[\frac{\alpha^{\prime}\beta}
{\alpha^{\prime}+p_i}+\frac{\alpha\beta^{\prime}}{\beta^{\prime}+p_j}\right],\
\ i,j\in{\mathbb N},\ \ {\rm where}\ \
\alpha,\alpha^{\prime},\beta,\beta^{\prime}>0,$$  if and only if $H$
is the  Freund BVE distribution with joint density
\begin{eqnarray*}
   { {h}(x,y)=\left\{\begin{array}{cc}
                    \alpha^{\prime}\beta e^{-(\alpha+\beta-\alpha^{\prime})y-\alpha^{\prime}x}, &
                    x>
                    y> 0\\\alpha\beta^{\prime}e^{-(\alpha+\beta-\beta^{\prime})x-\beta^{\prime}y}, & y\ge x> 0
                  \end{array}\right.}
\end{eqnarray*}
 (see Freund
1961).\medskip\\
{\bf Example 5.} The Laplace--Stieltjes transform ${\cal L}_H$
satisfies
$${\cal L}_H(p_i,p_j)=\frac{1}{(1+p_i)(1+p_j)-rp_ip_j},\ \
i,j\in{\mathbb N},\ \ {\rm where}\ \  r\in[0,1),$$   if and only if
$H$ is the Moran--Downton BVE distribution with density function
$$h(x,y)=\frac{1}{1-r}I_0\left(\frac{2\sqrt{rxy}}{1-r}\right)e^{-(x+y)/(1-r)},\ \ x,y>0,$$
in which $r$ is the correlation coefficient and
$I_0(t)=\sum_{k=0}^{\infty}(t/2)^{2k}/(k!)^2,\ t\in{\mathbb R}
:=(-\infty,\infty),$ is the modified Bessel function of the first
kind of order zero (see Moran 1967 and Downton 1970). Using ${\cal
L}_H,$ we see that the marginal distributions of $H$ have finite
second moments.
\medskip\\
 {\bf
Example 6.} The Laplace--Stieltjes transform ${\cal L}_H$ satisfies
$${\cal L}_H(p_i,p_j)=\frac{1}{\theta+p_i+p_j}\left[(\theta+p_i){\cal L}_F(p_i)+
(\theta+p_j){\cal L}_G(p_j)\right]-\frac{\theta}{\theta+p_i+p_j},\ \
i,j\in{\mathbb N},$$ where $\theta>0$ is a constant, if and only if
$H$ is the  bivariate lack-of-memory distribution $BLM(F,G,\theta)$
with survival function  (see Section 5 or Lin et al.\,2019):
\begin{eqnarray*}
   { \overline{H}(x,y)=\left\{\begin{array}{cc}
                    e^{-\theta y}\overline{F}(x-y), &\ x\ge y\ge 0 \\
                    e^{-\theta x}\overline{G}(y-x),  &\ y\ge x\ge 0.
                  \end{array}\right.}
\end{eqnarray*}
\noindent {\bf Example 7.} The Laplace--Stieltjes transform ${\cal
L}_H$ satisfies
$${\cal L}_H(p_i,p_j)=\frac{(1-r)^q}{(1-r+p_i+p_j+p_ip_j)^{q}},\ \
i,j\in{\mathbb N},\ \ {\rm where}\ \  r\in[0,1)\ {\rm and}\  q>0,$$  if and only if $H$
is the standard bivariate Gamma distribution with density function
$$h(x,y)=(1-r)^q\frac{(xy)^{q-1}}{\Gamma(q)}f_q(rxy)e^{-x-y},\ \ x,y>0,$$
in which $f_q(z)=\sum_{n=0}^{\infty}z^n/[n!\Gamma(q+n)],$ and $r,\
q>0$ are the correlation coefficient
 and the shape parameter, respectively (Letac et al.\,2007, p.\,14; Marcus 2014, Theorem 1.1).
\medskip\\
{\bf Example 8.} Let $H$ be the joint distribution of three
nonnegative variables $(X,Y,Z).$  The Laplace--Stieltjes transform
${\cal L}_H$ satisfies
\begin{eqnarray*}{\cal L}_H(p_i,p_j,p_k)&=&\bigg[\frac{1-a^2-b^2}{(1+p_i)(1+p_j)(1+p_k)}\bigg]^{\alpha}\\
&~~&\times \sum_{n
=0}^{\infty}\sum_{\ell=0}^{n}\frac{\Gamma(n+\alpha)}{\Gamma(\alpha)n!}{n\choose\ell}
\frac{a^{2\ell}b^{2(n-\ell)}}{(1+p_i)^{\ell}(1+p_j)^n(1+p_k)^{n-\ell}}
,\ \ i,j,k\in{\mathbb N},
\end{eqnarray*}
where $\alpha, a,b>0$ and $a^2+b^2<1$, if and only if $H$ is the
three-variate Gamma distribution with density function (see Marcus
2014, Example 2.3):
\begin{eqnarray*}
h(x,y,z)&=&(1-a^2-b^2)^{\alpha}(xyz)^{\alpha-1}e^{-(x+y+z)}\\
&~~&\times\sum_{n
=0}^{\infty}\sum_{\ell=0}^{n}\frac{y^n}{\Gamma(\alpha)n!}{n\choose\ell}
\frac{(a^2x)^{\ell}(b^2z)^{n-\ell}}{\Gamma(\ell+\alpha)\Gamma(n-\ell+\alpha)},\
\ x,y,z>0.
\end{eqnarray*}\\
\noindent {\bf Appendix}\\
{\bf Proof of the identity (16).} Rewrite the LHS of the
Lebesgue--Stieltjes integral in (15) as follows:
\begin{eqnarray}
&~&{\bf E}\bigg[\prod_{i=1}^{m+1} (1-\exp\big({-s_iX_i})\big)\bigg]\nonumber\\
&=&1+\sum_{k=1}^{m+1}(-1)^{k}{\bf E}\bigg[\sum_{1\le i_1<i_2<\cdots<i_k\le m+1}\exp(-s_{i_1}X_{i_1}-s_{i_2}X_{i_2}-\cdots-s_{i_k}X_{i_k})\bigg]\nonumber\\
&=&1+\sum_{k=1}^{m+1}(-1)^{k}\sum_{1\le i_1<i_2<\cdots<i_k\le m+1}{\bf E}[\exp(-s_{i_1}X_{i_1}-s_{i_2}X_{i_2}-\cdots-s_{i_k}X_{i_k})]\nonumber\\
&=&1+\sum_{k=1}^{m+1}(-1)^{k}\sum_{1\le i_1<i_2<\cdots<i_k\le m+1}{\cal L}_{H_{i_1i_2\ldots
i_k}}(s_{i_1},s_{i_2},\ldots,s_{i_k})\nonumber\\
&=&1+\sum_{k=1}^{m}(-1)^{k}\sum_{1\le i_1<i_2<\cdots<i_k\le m+1}\bigg(\prod_{j=1}^{k}s_{i_j}\bigg){L}_{H_{i_1i_2\ldots
i_k}}(s_{i_1},s_{i_2},\ldots,s_{i_k})\nonumber\\
&~&+(-1)^{m+1}{\cal L}_H(s_1,s_2,\ldots,s_{m+1}).
\end{eqnarray}
In the last equality, we use the assumption on induction: the
identity (11) holds true for $n=1,2,\ldots,m.$ Namely, for
$k=1,2,\ldots, m,$
\[{\cal L}_{H_{i_1i_2\ldots
i_k}}(s_{i_1},s_{i_2},\ldots,s_{i_k})=\bigg(\prod_{j=1}^{k}s_{i_j}\bigg){L}_{H_{i_1i_2\ldots
i_k}}(s_{i_1},s_{i_2},\ldots,s_{i_k}).
\]

 On the other hand, using (12) we rewrite the RHS
of the identity (15) as
\begin{eqnarray}
&~&\bigg(\prod_{i=1}^{m+1}s_i\bigg)\int_0^{\infty}\!\!\int_0^{\infty}\!\!\cdots\int_0^{\infty}
\overline{H}(x_1,x_2,\ldots,x_{m+1})\exp\Big({-\sum_{i=1}^{m+1}s_ix_i}\Big)\prod_{i=1}^{m+1}{\rm d}x_i\nonumber\\
&=&1+\sum_{k=1}^{m+1}(-1)^{k}\sum_{1\le i_1<i_2<\cdots<i_k\le m+1}\bigg(\prod_{i=1}^{m+1}s_i\bigg)\nonumber\\
&~~&\times\int_0^{\infty}\!\!\int_0^{\infty}\!\!\cdots\int_0^{\infty}H_{i_1i_2\ldots
i_k}(x_{i_1},x_{i_2},\ldots,x_{i_k})\exp\Big({-\sum_{i=1}^{m+1}s_ix_i}\Big)\prod_{i=1}^{m+1}{\rm d}x_i\\
&=&1+\sum_{k=1}^{m}(-1)^{k}\sum_{1\le i_1<i_2<\cdots<i_k\le m+1}\bigg(\prod_{j=1}^{k}s_{i_j}\bigg)L_{H_{i_1i_2\ldots
i_k}}(s_{i_1},s_{i_2},\ldots,s_{i_k})\nonumber\\
&~~&+(-1)^{m+1}\bigg(\prod_{i=1}^{m+1}s_i\bigg)L_H(s_1,s_2,\ldots,s_{m+1}).
\end{eqnarray}
In Eq.\,(18), we use the fact that
\[\bigg(\prod_{i=1}^{m+1}s_i\bigg)\int_0^{\infty}\!\!\int_0^{\infty}\!\!\cdots\int_0^{\infty}
\exp\Big({-\sum_{i=1}^{m+1}s_ix_i}\Big)\prod_{i=1}^{m+1}{\rm d}x_i=1,
\]
while in Eq.\,(19), we apply the result
\[\bigg(\prod_{i\in J_k}s_{i}\bigg)\int_0^{\infty}\!\!\int_0^{\infty}\!\!\cdots\int_0^{\infty}
\exp\Big({-\sum_{i\in J_k}s_ix_i}\Big)\prod_{i\in J_k}{\rm d}x_i=1,
\]
where $J_k=\{1,2,\ldots,m+1\}\setminus \{i_1,i_2,\ldots,i_k\},\
k=1,2,\ldots,m.$

Finally, comparing Eqs.\,(15), (17) and (19), we claim the identity
(11), or, equivalently, the identity (16). This completes the proof.
\medskip\\
\noindent {\bf Acknowledgments} The authors would like to thank the
Editor in Chief, Associate Editor and two referees for helpful
comments and suggestions, which improve the presentation of the
manuscript.
\bigskip\\
\noindent{\bf References}
\begin{description}
\item Aghili, A. and Parsa Moghaddam, B. (2011). Certain theorems on two
dimensional Laplace transform and non-homogeneous parabolic
partial differential equations. {\it Surv. Math. Appl.} {\bf 6},
165--174.

\item Almira, J.\,M. (2007). M\"untz type
theorems. I. {\it Surv. Approx. Theory} {\bf 3}, 152--194.

\item Apostol, T.\,M. (1975). {\it Mathematical Analysis}, 2nd edn.  Addison-Wesley,
Reading, MA.

\item Apostol, T.\,M. (1976). {\it Introduction to Analytic Number Theory}.  Springer, New York.

\item Balakrishnan, N. and Lai, C.-D. (2009). {\it Continuous Bivariate Distributions,} 2nd edn.
Springer, New York.

\item {Barlow, R.\,E. and Proschan, F. (1981)}. {\it
Statistical Theory of Reliability and Life Testing: Probability
Models}, To Begin With. Silver Spring, MD.


\item Boas, R.\,P. (1954). {\it Entire Functions.} Academic Press, New York.


\item Carson, J.\,R. (1919). Theory of the transient oscillations of
electrical networks and transmission systems. {\it Proceedings
of the American Institute of Electrical Engineers} {\bf 38},
407--489; {\it Transactions of the American Institute of
Electrical Engineers}, {\bf 38}, 345--427.

\item Cohen, A.\,M. (2007). {\it Numerical Methods for Laplace Transform Inversion:
Numerical Methods and Algorithms.} Springer, New York.

\item Ditkin, V.\,A. and Prudnikov, A. P. (1962). {\it Operational
Calculus in Two Variables and its Applications}, 2nd edn.
Academic Press, New York.

\item Donolato, C. (2002). Analytical and numerical inversion of the Laplace--Carson
transform by a differential method. {\it Comput. Phys. Commun.}
{\bf 145},  298--309.
\item Downton, F. (1970). Bivariate exponential distributions in reliability theory.
{\it J. Royal Statist. Soc. Ser. {B}}\ \ {\bf 32}, 408--417.

\item  Erd\'elyi, T. (2005).
The ``full M\"untz theorem" revisited. {\it Constr. Approx.}
{\bf 21}, 319--335.

\item  Erd\'elyi, T. and Johnson, W.\,B. (2001). The ``Full M\"untz Theorem" in $L_p[0, 1]$ for $0<p<\infty.$
{\it J. Anal. Math.} {\bf 84}, 145--172.

\item Farrell, R.\,H. (1976). {\it Techniques of Multivariate Calculations}. Springer, New York.
\item Feller, W. (1971). {\it An Introduction to Probability Theory and its Applications}, Vol. II, 2nd edn. Wiley, New York.
\item {Freund, J.\,E. (1961)}. A bivariate extension of the
exponential distribution. {\it J. Amer. Statist. Assoc.} {\bf
56}, 971--977.
\item Hougaard, P. (1986).
Survival models for heterogeneous populations derived from
stable distributions. {\it Biometrika} {\bf 73},  387-396.
\item Lerch, M. (1903). Sur un point de la th\'eorie de fonctions g\'en\'eratrices d'Abel.
{\it Acta Math.} {\bf 27}, 339--351.
\item Letac, G., Tourneret, J.-Y. and Chatelain, F. (2007).
The bivariate gamma distribution with Laplace transform $(1+as+bt+cst)^{-q}$:
 History, characterizations, estimation. Available online at
http://www.math.univ-toulouse.fr/$\sim$letac/bivgamma.pdf.
\item Lin, G.\,D. (1993). Characterizations of the exponential distribution via
the blocking time in a queueing system. {\it Statist. Sinica}
{\bf 3}, 577--581.
\item Lin, G.\,D. (1998). Characterizations of the ${\cal L}$-class of life distributions.
{\it Statist.} {\it Probab. Lett.} {\bf 40}, 259--266.
\item Lin, G.\,D., Dou, X. and Kuriki, S. (2019).
The bivariate lack-of-memory distributions. {\it Sankhy$\bar{a}$
Ser.\,A}\ \,{\bf 81}, 273--297.
\item  Lin, G.\,D., Lai, C.-D. and Govindaraju, K. (2016).
Correlation structure of the Marshall--Olkin bivariate
exponential distribution. {\it Statistical Methodology} {\bf
29}, 1--9.

\item Marcus, M.\,B. (2014). Multivariate gamma distributions. {\it Electron. Commun. Probab.} {\bf 19}, no. 86, 1--10.

\item {Marshall, A.\,W. and Olkin, I. (1967)}.
 A multivariate exponential distribution. {\em J. Amer. Statist.
 Assoc.} {\bf 62}, 30--44.
\item Miller, P.\,D. (2006). {\it Applied Asymptotic Analysis}.
American Mathematical Society,  Providence, RI.
\item Moran, P.\,A.\,P. (1967). Testing for correlation between
non-negative variates. {\it Biometrika} {\bf 54}, 385--394.
\item M\"untz, Ch.\,H. (1914). \"Uber den approximationssatz
von Weierstrass. In: H.\,A.\,Schwarz's Festschrift, Berlin,
pp.\,303--312.

\item P\'erez, D. and Quintana, Y. (2008). A survey on the Weierstrass approximation theorem.
{\it Divulg. Mat.} {\bf 16}, 231--247.

\item Ross, S.\,M. (2014). {\it Introduction to Probability Models.} Academic Press, Oxford.

\item  Rudin, W. (1987). {\it Real and Complex Analysis},  3rd edn.  McGraw-Hill, New York.

\item Sz\'asz, O. (1916). \"Uber die approximation stetiger funktionen durch
lineare aggregate von potenzen. {\it Math. Ann.} {\bf 77},
482--496.

\item van der Pol, B. and Bremmer, H. (1955). {\it Operational Calculus}, 2nd edn.
Cambridge University Press, Cambridge. Reprinted by Chelsea Press, New York (1987).

\item von Golitschek, M. (1983).
A short proof of M\"untz's theorem. {\it J. Approx. Theory} {\bf
39}, 394--395.

\item Weierstrass, K. (1885).  \"Uber die analytische
Darstellbarkeit sogenannter willk\"urlicher Functionen einer
reellen Ver\"anderlichen, {\it Verl. d. Kgl. Akad. d. Wiss.
Berlin} {\bf 2}, 633--639.

\item Young, W.\,H. (1917). On multiple integration by parts and the
second theorem of the mean. {\it Proc. Lond. Math. Soc.} {\bf
16}, 273--293.

\item Zaremba, S.\,K. (1968). Some applications of multidimensional integration by parts. {\it Ann. Pol. Math.} {\bf 21}, 85--96.

\end{description}
\end{document}